\input amstex
\frenchspacing
\documentstyle{amsppt}
%\magnification=\magstep1
\magnification=\magstephalf
\hsize=15.9 truecm \vsize=22.5 truecm
\footline{\hfill\fiverm version \today}
\def\today{\ifcase\month\or January \or February\or
March\or April\or May\or June\or July\or August\or
September\or October\or November\or December\fi
\space\number\day, \number\year}
\footline{\hfill\sevenrm version \today}
\NoBlackBoxes
\def\FF{\bold F}

\def\Aut{\mathop{\text{\rm Aut}}\nolimits}
\def\Gal{\mathop{\text{\rm Gal}}}
\def\GL{\mathop{\text{\rm GL}}\nolimits}
\def\disc{\mathop{\text{\rm disc}}}
\def\Hom{\mathop{\text{\rm Hom}}}

\def\isar{\ \smash{\mathop{\longrightarrow}\limits^\sim}\ }
\def\ab{{\text{\rm ab}}}
\def\tor{{\text{\rm tor}}}
\def\i{\mathop{\text{\rm i}}}
\def\ker{\mathop{\text{\rm ker}}}
\def\lcm{\mathop{\text{\rm lcm}}}
\def\ord{\mathop{\text{\rm ord}}}
\def\mapright#1{\ \smash{\mathop{\longrightarrow}\limits^{#1}}\ }
\def\legendre{\overwithdelims()}
\def\mod{\bmod}
\def\congr{ \equiv }
\def\iso{ \cong }
\def\tto{\longrightarrow}

\def\Que{\bold Q}

\def\Zee{\bold Z}

\def\muhat{\widehat{\mu}}
\newcount\refCount
\def\newref#1 {\advance\refCount by 1
\expandafter\edef\csname#1\endcsname{\the\refCount}}
\newref BR   %Brau
\newref CP   %Pappalardi
\newref HO   %Hooley
\newref JA   %Abtien scriptie
\newref JO   %Nathan Jones Serre curves
\newref LEa  %Lenstra 1977
\newref LEb  %Lenstra San Antonio notes
\newref MA   %Matthews
\newref MOa  %Moree 1998
\newref MOb  %Moree 2008
\newref MSa  %Moree-Stevenhagen 2-var Artin conjecture
\newref MSb  %Moree-Stevenhagen higher rank prim root densities
\newref WJP  %Palenstijn
\newref RZ   %Ribes-Zalesskii
\newref SCa  %Schinzel abelian binomials
\newref SCb  %Schinzel Selecta
\newref SE   %Serre open index
\newref ST   %psh on Artin
\newref WA   %Wagstaff
\topmatter
\title
Character sums for primitive root densities
\endtitle
\author H. W. Lenstra, Jr.,  P. Moree and  P. Stevenhagen\endauthor
\address 
Mathematisch Instituut, Universiteit Leiden, Postbus 9512, 2300 RA Leiden,
The Netherlands
\endaddress
\email hwl, psh\@math.leidenuniv.nl\endemail
\address 
Max-Planck-Institut f\"ur Mathematik,
Vivatsgasse 7, 53111 Bonn, Germany
\endaddress
\email moree\@mpim-bonn.mpg.de \endemail
%  Math Subject Classifications
\abstract
It follows from the work of Artin and Hooley that,
under assumption of the generalized Riemann hypothesis, the density
of the set of primes $q$ for which a given non-zero rational
number $r$ is a primitive root modulo $q$ can be written as 
an infinite product $\prod_p \delta_p$ of local factors $\delta_p$
reflecting the degree of the splitting field of $X^p-r$ at
the primes $p$, multiplied by a somewhat complicated factor
that corrects for the `entanglement' of these splitting fields.

We show how the correction factors arising in Artin's original primitive
root problem and several of its generalizations can be interpreted as
character sums describing the nature of the entanglement.
The resulting description in terms of local contributions
is so transparent that it greatly facilitates explicit 
computations, and naturally leads to non-vanishing criteria
for the correction factors.

The method not only applies in the setting of Galois representations
of the multiplicative group underlying Artin's conjecture, 
but also in the $\GL_2$-setting arising for elliptic curves.
As an application, we compute the density of the set of primes 
of cyclic reduction for Serre curves.
\endabstract
\endtopmatter

\document

\head 1. Introduction
\endhead

\noindent
\global\baselineskip=13pt
Artin's conjecture on primitive roots, which dates back to 1927,
states that for a non-zero rational number $r\ne \pm1$,
the set of primes $q$
with the property that $r$ is a primitive root modulo $q$
has natural density 
$$
\delta(r)=c_r\cdot\prod_{p\text{ prime}}
\Bigl(1-{1\over p(p-1)}\Bigr)
$$
inside the set of all primes.
Here $p$ ranges over the rational primes, and $c_r$ is a
rational number that depends on~$r$.
The universal constant
${\Cal A}=\prod_p (1-{1\over p(p-1)})\doteq .3739558$
is known as {\it Artin's constant}.
The constant $c_r$ vanishes if and only if $r$ is a square.
For values of $r$ that are not perfect powers, it was discovered
after computer calculations in 1957 that $c_r$ can be different from 1,
leading to a correction of the original conjecture by Artin himself [\ST].
In 1967, this corrected conjecture was proved under assumption of
the generalized Riemann hypothesis (GRH) by Hooley~[\HO].

The heuristic argument underlying Artin's conjecture is simple:
for a prime number $q$ that does not divide the numerator or
denominator of $r$, the number $r$ is a primitive root modulo $q$
if and only if there is no prime number $p$ dividing $q-1$ such that $r$
is a $p$-th power modulo $q$.
In terms of number fields, this condition amounts to saying
that for no prime $p<q$,
the prime $q$ splits completely in the splitting field
$$
F_p=\Que(\zeta_p, \root p\of r)\subset\overline\Que
\tag{1.1}
$$
of the polynomial $X^p-r$ over $\Que$.
Here $\overline\Que$ denotes an algebraic closure of $\Que$,
and~$\zeta_p$ a primitive $p$-th root of unity in $\overline\Que$.

For fixed~$p$, the set of primes $q$ that do not split completely
in $F_p$ has density $\delta_p=1-[F_p:\Que]^{-1}$.
As we have $r\ne\pm1$, there is a largest integer $h\in\Zee$ 
for which $r$ is an $h$-th power in $\Que^*$.
We have $[F_p:\Que]=p-1$ for $p$ dividing $h$,
and $[F_p:\Que]=p(p-1)$ otherwise.
If we assume that the splitting conditions at the various primes $p$
that we impose on $q$ are `independent', it seems reasonable to
conjecture, as Artin did, that $\delta(r)$ equals 
$$
A(r)=
\prod_p \delta_p=
%\prod_{p\text{ prime}} \left(1-{1\over [F_p:\Que]}\Bigr)=
\prod_p \Bigl(1-{1\over [F_p:\Que]}\Bigr)=
\prod_{p|h} \Bigl(1-{1\over p-1}\Bigr)\cdot
\prod_{p\nmid h} \Bigl(1-{1\over p(p-1)}\Bigr).\tag{1.2}
$$
Note that $A(r)$ is a rational multiple of Artin's constant $\Cal A$,
and equal to it for $h=1$.
We have $A(r)=0$ if and only if $r$ is a square; in this case
we have $\delta_2=0$, and $r$ is not a primitive root modulo any
odd prime $q$.

The assumption on the independence of the splitting
conditions in the various fields $F_p$ is not always correct.
If $F_2=\Que(\sqrt r)$ is a quadratic field of discriminant~$D$,
then it is contained in the $|D|$-th cyclotomic field $\Que(\zeta_{|D|})$. 
Thus, if $D$ is {\it odd\/}, then $F_2$ is contained in the compositum of
the fields $F_p$ with $p|D$.
This is however the only `entanglement' between the fields $F_p$
that occurs.
In order to take it into account,
one writes $F_n=\Que(\zeta_n, \root n\of r)$ for the splitting field
of $X^n-r$ and applies a standard inclusion-exclusion argument to obtain
a conjectural value
$$\delta(r)= \sum_{n=1}^\infty {\mu(n)\over [F_n:\Que]},\tag1.3$$
where $\mu$ denotes the M\"obius function.
If $F_2=\Que(\sqrt r)$ is not quadratic of odd discriminant, then
$[F_n:\Que]$ is a multiplicative function on squarefree values
of $n$, and (1.3) reduces to (1.2).
If $F_2$ is quadratic of odd discriminant $D$, then $[F_n:\Que]$
is no longer multiplicative, as it equals
${1\over 2}\prod_{p|n} [F_p:\Que]$ for all
squarefree~$n$ that are divisible by $D$.
In this case, it is not so clear whether the right hand side of
(1.3) is non-zero, or even non-negative.
However, a `rather harder' calculation [\HO, p.~219] shows that it
can be written as $\delta(r)=E(r) \cdot A(r)$, with $A(r)$ as in (1.2) and
$$
E(r)=1-\prod_{p|D\atop p|h} {-1\over p-2}\cdot
\prod_{p|D\atop p\nmid h} {-1\over p^2-p-1} \tag1.4
$$
an `entanglement correction factor'.
Note that $E(r)$ is well-defined as $D$ is odd.
The multiplicative structure of the second term of $E(r)$ immediately
shows that $E(r)$ is non-zero.

The explicit form of Artin's conjecture as we have just stated it,
is the version that was proved by Hooley
under the assumption of the generalized Riemann hypothesis.
The hypothesis is used to obtain sufficient control of the error
terms occurring in density statements for sets of primes that split completely
in the fields~$F_n$ in order to prove (1.3).
So far, unconditional results have remained insufficient to handle
conditions at infinitely many primes~$p$.

Artin's conjecture can be generalized in various ways.
For example, one may impose the additional condition on the primes $q$
that they lie in a given arithmetic progression.
Alternatively, one can replace the condition that $r$ be a primitive root
modulo~$q$ by the requirement that $r$ generate a subgroup of
given index in $\FF_q^*$, or even combine the two conditions.
Just as in the original conjecture, these conditions amount to imposing 
restrictions on the splitting behavior of $q$ in number fields
contained in the infinite extension
$$
F_\infty=\bigcup_{n=1}^\infty F_n=\bigcup_{n=1}^\infty
\Que(\zeta_n, \root n \of r)\subset\overline\Que
$$
of $\Que$.
They may be formulated as conditions on the Frobenius element
of $q$ in these number fields, for which density statements hold
by the Chebotarev density theorem.
As was shown by the first author [\LEa], the prime densities for such generalizations
can in principle (under assumption of GRH) be obtained
along the lines of Hooley's proof, and equal the `fraction' of
good Frobenius elements in $G$.
However, the explicit evaluation of the entanglement correction factor 
from a formula analogous to (1.3) rapidly becomes very unpleasant.

The present paper, which was already announced in [\ST], introduces a
simple but effective method to compute entanglement correction
factors for primitive root problems over $\Que$.
It starts with the observation 
that the Galois automorphisms of the field $F_\infty$ act as group
automorphisms on the subgroup
$$
R_\infty=\{x\in \overline\Que^*: x^k\in\langle r\rangle \subset\Que^*
          \text{ for some\ } k\in\Zee_{>0}\}
$$
of $\overline\Que^*$ consisting of the {\it radicals\/} 
that generate the infinite field extension $F_\infty$ of~$\Que$.

In Section~2, we show that for {\it all\/} $r\in\Que^*\setminus\{\pm1\}$,
this action gives rise to an
injective `Galois representation' 
$$
G=\Gal(F_\infty/\Que) \tto A= \Aut_{R_\infty\cap\Que^*} (R_\infty)
$$
that embeds $G$ as an open subgroup of index 2 in the group $A$ of
group automorphisms of $R_\infty$ that restrict to the identity
on $R_\infty\cap\Que^*$.
Unlike the Galois group $G$, the automorphism group $A$ is {\it always\/} a
product of local factors $A_p$ at the primes $p$.
In Theorem 2.9, we explicitly describe the quadratic character 
$\chi: A\to\{\pm1\}$ that has $G$ as its kernel: it is a finite
product $\chi=\prod_p\chi_p$ of quadratic characters $\chi_p$,
with each~$\chi_p$ factoring via the projection $A\to A_p$.

The profinite group $A=\prod_p A_p$ carries a Haar measure $\nu$,
which we can take equal to the product $\prod_p \nu_p$
of the normalized Haar measures $\nu_p$ on $A_p$.
For {\it any\/} subset $S\subset A$ of the form
$\prod_p S_p$ with $S_p\subset A_p$ measurable,
one can compute the fraction 
$\delta(S)=\nu(G\cap S)/\nu(G)$ of elements in $G$ that lie in $S$
as a character sum in terms of the quadratic character $\chi$.
In our applications, $S\cap G$ will be a set of `good' Frobenius
elements inside $G$. 
By Hooley's method, the fraction $\delta(S)$
is then, under GRH, the density of the set of
primes $q$ meeting the Frobenius conditions
imposed by the choice of $S$.

In Theorem 3.3, we show that for the sets $S=\prod_p S_p$ as above,
the fraction $\delta(S)$ is 
the natural product of a naive `Artin constant'
$$
\nu(S)=
{\nu(S)\over \nu(A)} = \prod_p {\nu_p(S_p)\over \nu_p(A_p)}
=\prod_p \nu_p(S_p)
$$
as we gave in (1.2) 
and an entanglement correction factor of the form
$$
E=1+\prod_p E_p.
\tag1.5
$$
Just as in (1.4), where we have $E_2=-1$, the
local factors $E_p$ are different
from 1 only for finitely many `critical' primes~$p$
occurring in the finite product $\chi=\prod_p\chi_p$.
The factor 
$$
E_p={1\over\nu_p(S_p)} \int_{S_p} \chi_p d\nu_p
$$
equals the average value of $\chi_p$ on $S_p$.
It is easily evaluated in cases where $S_p$
is a set-theoretic difference of subgroups of $A_p$,
and can usually be computed explicitly as the average value of
a quadratic character on a finite set.

The transparent structure of the formula obtained makes it easy
to decide when the fraction $\delta(S)$
of good Frobenius elements in $G$ vanishes.
Vanishing of the Artin constant $\nu(S)/\nu(A)$ means that
$S=\prod_p S_p$ is a set of measure 0.
In concrete examples, this implies that $S$ is empty, and that
there is a prime $p$ for which $S_p$ is empty because the conditions
imposed by $S$ cannot be met `at~$p$'.
In the original Artin conjecture, this only happens for $p=2$ in the
case that $r$ is a square.

Vanishing of the entanglement correction factor $E$
is a more subtle phenomenon that does not
occur in the original conjecture.
In accordance with Theorem~4.1 in~[\LEa], it means that there is an
incompatibility `at a finite level' between the conditions at the 
critical primes.
Since all $E_p$, being average values of characters, are bounded in
absolute value by 1, it is easy (cf. Corollary 3.4)
to spot the occurrences of $E=0$ in (1.5).
We illustrate this by computing the value $\delta(S)$ and its vanishing
criteria in the case of Artin's conjecture (Section 4) and
its generalizations to primes in arithmetic progressions (Section 5)
and near-primitive roots (Section 6) mentioned above.
These cases have been treated before, but the short statements and
complete proofs that we can now give both for the densities and
the vanishing criteria illustrates the strength and elegance of our method.

In Section 7 we briefly discuss extensions of our method that 
are currently being applied to more general Artin-like problems.
Including them with full details would have made this paper too long.
Instead, we prove a new result in Section 8 in a different direction,
in relation to so-called {\it Serre curves}.
These are elliptic curves over $\Que$ that are `generic' [\JO]
in the sense that their torsion points generate an infinite extension
of $\Que$ that is `as large as possible'.
For such curves, Serre proved that the density of primes $p$ for which 
the group $E(\FF_p)$ is cyclic is given (under GRH) by a formula
resembling (1.3).
We prove that this density can {\it also\/} be written as the product
of an elliptic Artin constant and a completely explicit
entanglement correction factor of the form (1.5).
It shows that the underlying idea of our method has a wide range
of application.

\head 2. The radical extension $F_\infty$
\endhead
\noindent
Let $r\in\Que^*$ be a non-zero rational number different from $\pm1$.
Then for $n\in\Zee_{\ge1}$, the number field
$F_n=\Que(R_n)=\Que(\zeta_n, \root n\of r)$ is obtained
by adjoining to $\Que$ a group of {\it radicals\/}
$$
R_n= \{x\in \overline\Que^*: x^n\in\langle r\rangle \subset\Que^*\}
\subset\overline\Que^*.
$$
As $R_n$ is stable under the action of Galois,
the action of a {\it field\/} automorphism on $F_n$ is completely
determined by its action as a
{\it group\/} automorphism on the group of radicals $R_n$.
This gives rise to a natural injection of finite groups
$$
\Gal(F_n/\Que) \tto A(n)= \Aut_{R_n\cap\Que^*} (R_n).
\tag{2.1}
$$
The union $R_\infty=\bigcup_{n\ge1} R_n$ of all radical groups
generates an infinite algebraic extension $F_\infty=\Que(R_\infty)$
of $\Que$, and we
may take projective limits on both sides of (2.1) to obtain
the injective group homomorphism
$$
G=\Gal(F_\infty/\Que) \tto A= \Aut_{R_\infty\cap\Que^*} (R_\infty)
\tag{2.2}
$$
announced in the Introduction.
Note that the profinite groups $G$ and $A$ each come equipped with
their Krull topology, and that the injection (2.2) is an injection of
{\it topological\/} groups.

As $R_1=\langle r\rangle$ is free of rank 1, we have
$$
R_\infty\cap\Que^* =\langle r_0\rangle \times \langle-1\rangle
\tag{2.3}
$$
for a rational number $r_0\ne 0, \pm 1$ that is uniquely determined up
to sign, and up to inversion.
We fix the exponent by writing
$r=\pm r_0^e$, with $e>0$ the index of the subgroup
$\langle r  \rangle \times \langle-1\rangle$
inside  $\langle r_0\rangle \times \langle-1\rangle$.
If $e$ is odd, or $r$ is positive,
it is possible to choose the sign of $r_0$ in such a way
that we have $r=r_0^e$.
In the {\it twisted case\/} however, where $-r$ is a square and
$\Que(\sqrt r)=\Que(\i)$ the Gaussian number field,
$r$~itself is not an $e$-th power, and we have $r=-r_0^e$.
As we will see in the Sections 4--6, this twisted case often needs 
special attention in explicit computations.

The group $R_\infty$ contains the group
$\mu_\infty=\bigcup_{n\ge1}\mu_n(\overline\Que)$
of roots of unity as a subgroup.
As $\overline\Que^*$ is a divisible group, we can extend
the group embedding $\Zee\isar \langle r_0 \rangle \subset \overline\Que^*$
defined by $1\mapsto r_0$
to an embedding $\Que \tto \overline\Que^*$.
Note that this amounts to giving a section for the natural map
$R_\infty \to R_\infty/\mu_\infty\iso\Que$.
Given~$r_0$, we fix such an embedding~$q\mapsto{r_0^q}$. 
We can then write $r_0^\Que$ for its image in $\overline\Que^*$,
and $r_0^{1/n}$ for the $n$-th root of $r_0$ in this image.
With this notation, we have
$$
R_\infty= r_0^\Que \times \mu_\infty.
\tag{2.4}
$$
The automorphism group $A= \Aut_{R_\infty\cap\Que^*} (R_\infty)$
comes with a natural restriction map $A\to\Aut(\mu_\infty)$
that is continuous and surjective, and that admits a continuous
left inverse: extend the action to be the identity on $r_0^\Que$.
As an automorphism
$\sigma$ of $R_\infty$ that is the identity on $\mu_\infty$ is determined by
the values $\sigma(r_0^{1/n})/r_0^{1/n}\in\mu_n(\overline\Que)$
for $n\ge1$, we deduce from (2.3) that the lower row in the commutative
diagram of topological groups below is a split exact sequence
describing $A$:
$$

\def\mapdown#1{\Big\downarrow\rlap{$\vcenter{\hbox{$\scriptstyle#1$}}$}}
\matrix
1&\mapright{}&\Gal(F_\infty/\Que_\ab)&\mapright{}&G
             &\mapright{}&\Gal(\Que_\ab/\Que)&\mapright{}&1\cr
        &&\mapdown{}&&\mapdown{(2.2)}&&\mapdown{\wr}&&\cr
        1&\mapright{}&\Hom(r_0^\Que/r_0^\Zee,\mu_\infty)&\mapright{}&A
             &\mapright{}&\Aut(\mu_\infty)&\mapright{}&1\rlap{.}\cr
\endmatrix
$$
\smallskip\noindent
The upper row is exact by Galois theory, and the right vertical
isomorphism reflects the fact
that the maximal cyclotomic extension $\Que(\mu_\infty)$ of $\Que$
is the maximal abelian extension $\Que_\ab$ of $\Que$, which has
Galois group $\Gal(\Que_\ab/\Que)=\Aut(\mu_\infty)$.
As all automorphisms of $\mu_\infty$ are exponentiations,
$\Aut(\mu_\infty)$ is isomorphic to
the unit group $\widehat\Zee^*$ of the profinite completion
$\widehat\Zee=\lim\limits_{\leftarrow n} (\Zee/n\Zee)$
of $\Zee$.

We see that, in order to understand $G$ as a subgroup of $A$, we need
to identify the image of the Kummer map
$$
\eqalign{
\Gal(F_\infty/\Que_\ab)&\tto \Hom(r_0^\Que/r_0^\Zee,\mu_\infty) \cr
\sigma       &\longmapsto [r_0^x\mapsto (r_0^x)^{\sigma-1} ]. \cr
}
$$
By Kummer theory, this image is the abelian group dual to
$r_0^\Que/(r_0^\Que\cap \Que_\ab^*)$.
\proclaim
{2.5. Lemma}
We have $r_0^\Que\cap \Que_\ab^*= r_0^{{1\over 2}\Zee}$,
%=\langle\sqrt{r_0}\rangle$,
and $G$ is an open subgroup of $A$ of index~$2$.
\endproclaim
\noindent
{\bf Proof.}
The equality for $r_0^\Que\cap \Que_\ab^*$ amounts to saying that the
largest integer $k$ for
which the splitting field of $X^k-r_0$ is abelian over $\Que$ equals 2.
As neither $r_0$ nor $-r_0$ is a perfect power in $\Que^*$,
this is an immediate corollary of Schinzel's theorem on abelian
binomials [\SCa, Theorem 2; \SCb],
which implies that this splitting field is abelian if and only if
$r_0^{\#\mu_k(\Que)}=r_0^2$ is a $k$-th power in $\Que^*$.

The group $r_0^\Que/r_0^{{1\over 2}\Zee}$ is the quotient of
$r_0^\Que/r_0^\Zee$ that is obtained by dividing out its unique
subgroup of order 2, which is generated by $r_0^{1/2}\mod r_0^\Zee$. 
Its $\mu_\infty$-dual is the open subgroup 
$\Hom(r_0^\Que/r_0^{(1/2)\Zee},\mu_\infty)$ of index~2
in $\Hom(r_0^\Que/r_0^\Zee,\mu_\infty)$, and
we conclude that $G$ is open and of index 2 in $A$.
\hfill$\square$
\medskip\noindent
Lemma 2.5 actually yields an explicit description of the subgroup $G\subset A$
that arises as the Galois group of $F_\infty$ over $\Que$.
The group
$\Hom(r_0^\Que/r_0^\Zee,\mu_\infty)\iso\Hom(\Que/\Zee, \mu_\infty)$
can be viewed as the Tate module 
$$
\muhat=\lim_{\leftarrow n}\mu_n=
\Hom(\lim_{\to n} (\textstyle{1\over n}\Zee/\Zee), \mu_\infty)
=\Hom(\Que/\Zee, \mu_\infty)
\tag{2.6}
$$
of the multiplicative group.
It is a free module of rank 1 over
$\widehat\Zee$, 
and the natural action of $\Aut(\mu_\infty)=\widehat\Zee^*$
on $\Hom(\Que/\Zee,\mu_\infty)$ via the second argument
is simply the $\widehat\Zee^*$-multiplication we have on the 
$\widehat\Zee$-module $\muhat$.

From the split exact sequence for $A$, we see
that $A$ is the semidirect product
$$
A=\Hom(r_0^\Que/r_0^\Zee,\mu_\infty)\rtimes \Aut(\mu_\infty)
=\muhat\rtimes \widehat\Zee^*,
\tag{2.7}
$$
which is isomorphic to the affine group
$\widehat\Zee\rtimes \widehat\Zee^*$ over $\widehat\Zee$.

The subgroup $G\subset A$ is an extension of $\widehat\Zee^*$
by $\muhat^2\iso\muhat$, but this extension is non-split: if it were,
the closure of the commutator subgroup of $G=\Gal(F_\infty/\Que)$
would be of index 2 in $\Gal(F_\infty/\Que_\ab)$, contradicting the
fact that $\Que_\ab$ is the maximal abelian extension of $\Que$.

By (2.4), the field $F_\infty=\Que(R_\infty)$ is the compositum of 
$\Que(r_0^\Que)$ and $\Que(\mu_\infty)=\Que_\ab$,
and the embedding $G\subset A$, with $A$ in the explicit form (2.7),
amounts to a description of the field automorphisms of $F_\infty$
in terms of their action on these constituents.
The index 2 of $G$ in $A$ reflects the fact that by Lemma 2.5,
the intersection of
$\Que(r_0^\Que)$ and $\Que_\ab$ is not $\Que(r_0)=\Que$,
but the quadratic field $K=\Que(\sqrt {r_0})$. 
This implies that an element $(\phi,\sigma)\in A$ is in $G$ if and only if 
$\phi\in \Hom(r_0^\Que/r_0^\Zee,\mu_\infty)$ and
$\sigma\in \Aut(\mu_\infty)=\Gal(\Que_\ab/\Que)$
act in a compatible way on $\sqrt{r_0}=r_0^{1/2}\in\Que_\ab$:
$$
\phi(r_0^{1/2}) = (r_0^{1/2})^{\sigma-1} 
\in\mu_2.
$$
In words: an automorphism of the multiplicative group of radicals $R_\infty$ 
induces an automorphism of the field $F_\infty=\Que(R_\infty)$
if and only if it respects the {\it additive\/} identity expressing
$\sqrt{r_0}=r_0^{1/2}$
as a $\Que$-linear combination of roots of unity with
rational coefficients.
We can phrase this slightly more formally by saying that $G\subset A$
is the difference kernel of two distinct quadratic characters 
$\psi_K, \chi_K: A\to\mu_2$ related to the
{\it entanglement field\/} $K=\Que(\sqrt{r_0})$.

The first quadratic character $\psi_K:A\to\{\pm 1\}$
describes the action on the generator $r_0^{1/2}$ of $K$
in terms of the $\phi$-component of $a=(\phi,\sigma)\in A$:
$$
\psi_K(a)=\phi(r_0^{1/2})\in\mu_2.
\tag{2.8}
$$
Note that $\psi_K$ is indeed a character on $A$, as $\Aut(\mu_\infty)$
acts trivially on $\mu_2=\{\pm1\}$.
It is a {\it non-cyclotomic\/} character on $A$, i.e., it does
not factor via the natural map
$\omega: A\to\Aut(\mu_\infty)$.
The second character
$$
\chi_K = {r_0 \overwithdelims() \cdot}:
A\mapright{\omega} \Aut(\mu_\infty)=\widehat\Zee^*\to \mu_2
$$
is a cyclotomic character on $A$, of conductor $d=\disc(K)$,
associated to the entanglement field $K=\Que(\sqrt{r_0})$.
It factors via the quotient $(\Zee/d\Zee)^*$ of the
cyclotomic component $\Aut(\mu_\infty)=\widehat\Zee^*$ of $A$,
on which it can be viewed as a Dirichlet character.
Its value in $a\in (\Zee/d\Zee)^*$ is given by the Kronecker symbol
${r_0 \overwithdelims() a}$ corresponding to~$K$.
\proclaim
{2.9. Theorem}
Let $K=\Que(\sqrt{r_0})$ and $\psi_K, \chi_K: A\to\mu_2$
be defined as above.
Then the natural map $(2.2)$ identifies
$G=\Gal(F_\infty/\Que)$ with the subgroup of $A$ of index $2$
that arises as the kernel of the quadratic character
$$
\psi_K\cdot \chi_K : A\to \mu_2.
\eqno{\square}
$$
\endproclaim
\noindent
From the description in (2.7) or, more canonically,
from the fact that automorphisms
of $R_\infty$ over $R_\infty\cap\Que^*$ can be given in terms of
their action on prime power radicals, it is clear that $A$ admits
a natural decomposition
$$
A=\prod_{p\text{ prime}} A_p, 
\tag{2.10}
$$
with $A_p$ the group of automorphisms of the group
$R_{p^\infty}=\bigcup_{k\ge1} R_{p^k}$
of $p$-power radicals that restrict to the identity on
$R_{p^\infty}\cap\Que^*$.

The character $\psi_K$ in Theorem 2.9 factors via the component
$A_2$ of $A$.
The other character $\chi_K$ can be decomposed
in the standard way for Dirichlet characters as a product
$$
\chi_K = \prod_p \chi_{K, p}
$$
of quadratic characters
$$
\chi_{K, p}: A\to \Aut(\mu_\infty)=\widehat\Zee^*\to\Zee_p^*\to \mu_2
\tag{2.11}
$$
of $p$-power conductor that are non-trivial exactly for primes $p$ 
dividing $d=\disc(K)$.
For odd primes $p|d$, the character $\chi_{K, p}$ is a lift to $A$ of
the Legendre symbol at $p$. 
For $d$ even, the character $\chi_{K, 2}$ is a lift to $A$
of one of the three non-trivial Dirichlet characters of conductor
dividing~8.
Note that $\chi_{K, p}$ factors via $A_p$ for all $p$.
\medskip\noindent
{\bf 2.12. Remarks.}
{\bf 1.}
As the rational number $r_0$ in (2.3) is only determined up to sign, the
entanglement field $K=\Que(\sqrt {r_0})$ in Theorem 2.9 
is only unique up to twisting by the cyclotomic character
$\chi_{-4}: A \to \Aut(\mu_\infty)=\widehat\Zee^*\to\{\pm 1\}$
giving the action on $\i=\sqrt{-1}$.
In other words, $G$ is {\it also\/} the subgroup of $A$ on which the
quadratic characters 
$\psi_{K'}=\psi_K\cdot\chi_{-4}$ and 
$\chi_{K'}=\chi_K\cdot\chi_{-4}$ 
related to $K'=\Que(\sqrt{-r_0})$ coincide.
Indeed, as $\chi_{-4}$ is quadratic,
the two products $\psi_{K'}\cdot \chi_{K'}$ and $\psi_K\cdot \chi_K$
defining~$G$ are equal.

{\bf 2.}
For odd $p$, one may identify the $p$-component $A_p$ of $A$
with the Galois group of $F_{p^\infty}=\Que(R_{p^\infty})$ over $\Que$.
For $p=2$, this is only true if we are not in the special case where
the entanglement field $K=\Que(\sqrt{r_0})$ equals $\Que(\sqrt {\pm2})$, and
$\psi_K\cdot \chi_K$ factors via $A_2$.
In non-special cases, there is a true `entanglement' of the
extensions $F_{p^\infty}$ in the sense that the character $\psi_K\cdot\chi_K$
in Theorem 2.9 that determines $G$ as a subgroup of $A$
is non-trivial on more than one prime component~$A_p$.
In the special case $K=\Que(\sqrt {\pm2})$, we have $d=\pm 8$ and
there is no entanglement at the level of Galois groups;
we do however have $G=G_2\times\prod_{p>2} A_p$ for a subgroup
$G_2\subset A_2$ of index 2.

{\bf 3.}
As we saw in (2.7), the semidirect product $A=\muhat\rtimes\widehat\Zee^*$
is a split extension of $\widehat\Zee^*$ by the multiplicative
Tate module $\muhat$ from (2.6).
The subgroup $G\subset A$ `cut out' by $r_0$ in the sense of 
Theorem 2.9 is not.
It is a {\it non-split\/} extension of $\widehat\Zee^*$ by 
the subgroup $\muhat^2\subset \muhat$ of index 2, which
is again isomorphic to $\muhat$.
Even though this is not directly relevant for us, one may wonder
which non-split extensions 
$$
\varepsilon_{r_0}: \qquad
1\to\muhat\tto G\tto
\widehat\Zee^*\to 1
$$
of profinite groups 
are provided by the Galois groups $G=\Gal(F_\infty/\Que)$
for various choices of $r_0\in\Que^*/\{\pm1\}$.
The answer is that {\it every\/} non-split extension
of $\widehat\Zee^*$ by $\muhat$ arises in this way, for a 
quadratic field $\Que(\sqrt {r_0})$ that is unique up to twisting
by $\chi{-4}$ in the sense of Remark 2.12.1.
A more formal way to phrase this consists in the construction of
an isomorphism
$$
\Que^*/(\{\pm1\}\cdot(\Que^*)^2) \isar H^2(\widehat\Zee^*,\muhat)
$$
under which $\pm r_0 {\Que^*}^2$ corresponds to the class
of the extension $\varepsilon_{r_0}$ in a
continuous cochain cohomology group $H^2(\widehat\Zee^*,\muhat)$
that describes profinite group extensions 
in the spirit of [\RZ, Theorem 6.8.4].
Such a construction can be given by standard arguments considering the
$\widehat\Zee^*$-cohomology of the sequence
$1\to\muhat\mapright{\square}\muhat\tto\mu_2\to 1$ describing multiplication
by~2 on $\muhat$, but one first has to establish the necessary formal
properties of continuous cochain cohomology groups $H^q(G, A)$
for profinite rather than simply discrete $G$-modules~$A$.
This is achieved in the Leiden master thesis [\JA] of Abtien Javanpeykar.

\bigskip

\head 3. Entanglement correction using character sums
\endhead

\noindent
The automorphism group $A$ and each of its components $A_p$ in (2.10) are
infinite profinite groups that naturally come with a topology
and a Haar measure.
The quadratic character $\psi_K\cdot\chi_K$ in Theorem 2.9 is continuous on $A$
with respect to this topology, and $G$ is an open subgroup of $A$
of index 2.
We normalize the Haar measure~$\nu_p$ on the compact groups~$A_p$
by putting $\nu_p(A_p)=1$;
this makes the product measure $\nu=\prod_p \nu_p$ into
a normalized Haar measure on $A$.

Densities for Artin-like primitive root problems (in one generator
over $\Que$) arise as fractions
$\delta(S)=\nu(G\cap S)/\nu(G)$ of `good' Frobenius elements
inside the Galois group $G=\Gal(F_\infty/\Que)$ of Theorem 2.9.
Here 
$$
S=\prod_p S_p\subset \prod_p A_p=A
$$
is some measurable subset of $A$ that
is defined componentwise at each prime $p$. 
Usually $S_p$ is the inverse image of some finite set $\overline S_p\subset 
\overline A_p$ under a continuous map $A_p\to \overline A_p$
to a finite discrete group $\overline A_p$.
A frequently encountered example is, for $P$ a power of $p$,
the restriction map
$$
\varphi_P: A_p \tto A(P) =\Aut_{R_P\cap\Que^*}(R_P).
\tag{3.1}
$$
Note that unlike $R_\infty$, the group $R_P$ of all $P$-th roots
of $\langle r\rangle$ depends on $r$, not just on $r_0$.
For Artin's original conjecture, the condition at~$p$
on the Frobenius element is that it is non-trivial on the field
$F_p=\Que(R_p)=\Que(\zeta_p, \root p \of r)$ from (1.1),
so we take
$$
\varphi_p: A_p \tto A(p) =\Aut_{R_p\cap\Que^*}(R_p)
%\iso\Gal(F_p/\Que)
$$
with $\overline S_p= A(p)\setminus\{1\}$, and put
$$
S_p=A_p\setminus \ker \varphi_p = \varphi^{-1}[\overline S_p].
\tag{3.2}
$$
As $\varphi_p$ is surjective and $A(p)\iso\Gal(F_p/\Que)$
has order $[F_p:\Que]$,
the subset $S_p\subset A_p$ has measure
$\nu_p(S_p)= 1 - [F_p:\Que]^{-1}$.
Thus, in this case $S=\prod_p S_p$ has
measure
$$
\nu(S)=\prod_p \nu_p(S_p)=
\prod_p \Bigl(1-{1\over [F_p:\Que]}\Bigr)
$$
equal to the constant $A(r)$ occurring in (1.2).
The entanglement correction factor $E(r)$ in (1.4)
is the factor by which 
$\delta(S)=\nu(G\cap S)/\nu(G)$ is different from $\nu(S)=\nu(S)/\nu(A)$
for the subgroup $G\subset A$ of index 2 described by Theorem~2.9.
Such entanglement correction factors can be computed in
great generality from the following theorem.
\proclaim
{3.3. Theorem}
Let $A=\prod_p A_p$ be as in $(2.10)$, with 
Haar measure $\nu=\prod_p \nu_p$, and 
$\chi=\prod_p \chi_p: A\to\{\pm1\}$ a non-trivial character obtained
from a family of continuous quadratic characters
$\chi_p: A_p\to \{\pm1\}$, with $\chi_p$ trivial for almost all
primes $p$.
Then for $G=\ker\chi$ and
$S=\prod_p S_p$ a product of $\nu_p$-measurable subsets 
$S_p\subset A_p$ with $\nu_p(S_p)>0$, we have 
$$
\delta(S)=
{\nu(G\cap S)\over \nu(G)} = 
\Bigl(1+\prod_p E_p\Bigr) \cdot {\nu(S)\over \nu(A)},
$$
with $E_p= E_p(S)={1\over \nu_p(S_p)}\int_{S_p}\chi_p d\nu_p$
the average value of $\chi_p$ on $S_p$.
\endproclaim\noindent
{\bf Proof.}
We assume that $\nu(S)=\prod_p \nu_p(S_p)$ is positive, as the theorem
trivially holds for $\nu(S)=0$.
We compute $\nu(G\cap S)$ by integrating the characteristic function
$(1+\chi)/2$ of $G$ over the subset $S\subset A$ with respect to $\nu$.
As we have $\nu(G)={1\over 2}\nu(A)$ by the non-triviality of $\chi$,
we obtain
$$
{\nu(G\cap S)\over \nu(G)} =
{1\over\nu(A)} \int_S (1+\chi) d\nu =
{\nu(S)\over \nu(A)}\cdot\left(1+{1\over\nu(S)}\int_S \chi d\nu\right).
$$
Now $\nu(S)$ equals $\prod_p \nu_p(S_p)$,
and the integral of $\chi=\prod_p \chi_p$ over $S=\prod_p S_p$ is the product
of the values $\int_{S_p}\chi_p d\nu_p$ for all $p$.
\qed
\proclaim
{3.4. Corollary}
For a set $S$ of positive measure,
the density $\delta(S)$ in $3.3$ vanishes 
if and only if 
there exists a sequence 
$\{\varepsilon_p\}_p$ of signs $\varepsilon_p\in\{\pm1\}$,
almost all equal to~$1$, such that we have
$\prod_p \varepsilon_p =-1$,
and
$\chi_p=\varepsilon_p$ almost everywhere on $S_p$.
\endproclaim\noindent
{\bf Proof.}
Suppose we have $\delta(S)=0$ and $\nu(S)>0$.
This amounts to saying that
the product $\prod_p E_p$, which is finite as we have 
$E_p=1$ for all $p$ at which $\chi_p$ is trivial, equals $-1$.
As every $E_p$ is the average value of a quadratic character on $S_p$,
it is a real number in $[-1, 1]$. 
It equals 1 (or $-1$) if and only if $\chi_p$ is
$\nu_p$-almost everywhere equal to 1 (or $-1$) on~$S_p$.
Thus, $\prod_p E_p =-1$ occurs if and only if the conditions
listed are satisfied.
\qed
\medskip\noindent
For the Galois group $G=\Gal(F_\infty/\Que)\subset A$ from
Theorem 2.9, we are in the situation of Theorem 3.3 in view of (2.11):
take $\chi=\prod_p\chi_p=\psi_K\cdot\chi_K$ with
$$
\chi_p=\cases
\psi_K\cdot\chi_{K,2}&\text{for $p=2$;}\cr
\chi_{K,p}&\text{for $p>2$.}\cr
\endcases
\tag{3.5}
$$
The characters $\chi_{K,2}$ and $\psi_K$ cannot coincide on $A$, as 
$\chi_{K,2}$ factors via the cyclotomic component $\widehat\Zee^*$
of $A$ in (2.7), whereas $\psi_K$ does not.
It follows that $\chi_2$ is always non-trivial.
Note also that, just as in Remark 2.12.1, the character
$\chi_2$ is unchanged if we replace 
$K=\Que(\sqrt{r_0})$ by $K'=\Que(\sqrt{-r_0})$.
\medskip\noindent
{\bf 3.6. Remark.}
As we noticed in Remark 2.12.2, it can happen in the situation of
Theorem 3.3 that all $\chi_p$'s but one character $\chi_q$ are trivial.
In this case we have $G=G_q\times \prod_{p\ne q} A_p$ for some
subgroup $G_q\subset A_q$ of index 2, and
$G\cap S$ will be the same for all subsets $S_q\subset A_q$
having the same intersection $S_q'=S_q\cap G_q$.
The correction factor $1+\prod_p E_p=1+E_q$ does however depend on
$S_q$, not only on $S_q'$.
This is not a contradiction, since
we can write $S_q=S_q'\cup S_q''$ as a disjoint union
with $S_q''=S_q\cap (A_q\setminus G_q)$,
and observe that the right hand side in Theorem 3.3 equals
$$
\eqalignno{
\Big(1+\prod_p E_p\Big){\nu(S)\over\nu(A)}&=
{1\over\nu(A)}
\Big(\nu_q(S_q)+\int_{S_q}\chi_q d\nu_q\Big)
\prod_{p\ne q}\nu_p(S_p)\cr
&={1\over\nu(A)}
\left(\nu_q(S_q)+\nu_q(S_q')-\nu_q(S_q'') \right)
\prod_{p\ne q}\nu_p(S_p)\cr
&={1\over\nu(G)}\nu_q(S_q')\prod_{p\ne q}\nu_p(S_p),\cr}
$$
in accordance with the fact that we have
$G\cap S=S_q'\times\prod_{p\ne q}S_p$.

\head 4. Artin's conjecture
\endhead

\noindent
Theorems 2.9 and 3.3 reduce the computation of the correction
factors occurring in many Artin-like problems
to fairly mechanical computations.
For Artin's original problem, which only takes a rational number
$r\in\Que^*\setminus\{\pm1\}$ as its input, we already noticed in (3.2)
that each subset $S_p\subset A_p$ of `good' Frobenius elements 
at $p$ equals $S_p= A_p\setminus\ker\varphi_p$ for
the natural map 
$$
\varphi_p: A_p \to A(p) =\Aut_{R_p\cap\Que^*}(R_p)\iso\Gal(F_p/\Que).
$$
This gives rise to the {\it Artin set\/} $S=S(r)=\prod_p S_p$, which
has (normalized) measure $\nu(S)=A(r)$ inside $A$ given by (1.2).
We have $\nu(S)=0$ if and only if $r$ is a square in $\Que^*$;
in this case, $S$ is empty as we have $S_2=\emptyset$.

To recover the correction factor $E(r)$ from (1.4) for non-square $r$,
we need to compute the entanglement correction factor $1+\prod_p E_p$
established in Theorem 3.3.
As $S_p=A_p\setminus \ker\varphi_p$ is the set-theoretic difference
of a group and a subgroup, the average value
$$
E_p=
{1\over\nu(S_p)} \left[\int_{A_p} \chi_p\,d\nu_p-
                       \int_{\ker\varphi_p} \chi_p\,d\nu_p\right]
$$
of $\chi_p$ on $S_p$
can only have three possible values, depending on the nature of $\chi_p$.
If $\chi_p$ is trivial, we obviously have $E_p=1$.
If $\chi_p$ is non-trivial on $\ker\varphi_p$, and therefore on $S_p$,
we get $E_p=0$ as both integrals, being integrals of
a non-trivial character over a group, vanish.
The interesting case is where $\chi_p$ is trivial on $\ker\varphi_p$
but not on $A_p$, and $E_p$ assumes the value
$$
E_p={-\nu_p(\ker\varphi_p)\over\nu_p(S_p)}
    ={-[F_p:\Que]^{-1}\over 1-[F_p:\Que]^{-1}}
    ={-1\over [F_p:\Que]-1}.
\tag{4.1}
$$
As we have $[F_p:\Que]=p-1$ if $r$ is a $p$-th power in $\Que^*$,
and $[F_p:\Que]=p^2-p$ otherwise, these $E_p$ are exactly the factors
that occur, for $D=\disc(F_2)$ odd and $p|2D$,
in the correction factor $E(r)$ in~(1.4).
The value of $E_2\in\{0, -1\}$ actually depends on the parity of $D$,
and the density correction for Artin's primitive root conjecture
can be formulated as follows.
\proclaim
{4.2. Theorem}
Let $r\ne-1$ be a non-square rational number, 
$G\subset A$ as in Theorem~$2.9$, and $S=S(r)\subset A$
the Artin set defined above.
Then $S$ has measure $A(r)$ given by $(1.2)$, and we have 
$$
\delta(S)={\nu(G\cap S)\over \nu(G)}= E(r)\cdot A(r)
$$
for an entanglement correction factor $E(r)$ that has the value $1$
if $D=\disc(\Que(\sqrt r))$ is even, and 
the value
$$
E(r)=1 +\prod_{p|2D} {-1\over [F_p:\Que]-1}
$$
from $(1.4)$ if $D$ is odd.
\endproclaim
\noindent
{\bf Proof.}
We apply Theorem 3.3, with $\chi=\psi_K\cdot\chi_K$ the character
from Theorem 2.9. Here we have $K=\Que(\sqrt{r_0})$, with $r=\pm r_0^e$
defined as in (2.3).
As we know already that $S$ has measure $A(r)$, we only have to
compute the factors $E_p$ occurring in the correction
factor $E(r)=1+\prod_p E_p$.
In our case, $E_p$ is the average value of the character $\chi_p$
from (3.5) on the set $S_p=A_p\setminus \ker\varphi_p$ from (3.2).

Suppose first that we are {\it not\/} in the twisted case where $-r$
is a square.
Then we can take $r=r_0^e$ with $e$ odd, and the fields 
$F_2=\Que(\sqrt r)$ and $K=\Que(\sqrt{r_0})$ coincide.
The character $\chi=\prod_p \chi_p$ has non-trivial
$p$-components only at $p=2$ and at the odd primes $p$ dividing
$D=\disc (K)$.
At odd primes $p|D$, the Legendre symbol $\chi_p$
is trivial on $\ker\varphi_p$ but not on $A_p$, so $E_p$ is given by (4.1).
At $p=2$, the non-trivial character $\chi_2=\psi_K\cdot\chi_{K,2}$ on $A_2$
equals $\chi_{K,2}$ on $\ker\varphi_2=\ker\psi_K$.
If $D$ is odd, it is trivial on $\ker\varphi_2$ and we find
$E_2=-1/([F_2:\Que]-1)=-1$ from (4.1), yielding $E(r)$ as stated.
If $D$ is even, $\chi_{K, 2}$ and therefore $\chi_2$ are
non-trivial on $\ker\varphi_2$, since $\varphi_2$ is non-cyclotomic;
we find $E_2=0$ and $E(r)=1$.

In the twisted case $D=-4$ the field $F_2=\Que(\sqrt r)=\Que(\i)$
is different from $K=\Que(\sqrt{r_0})$, and the character
$\varphi_2=\chi_{-4}$ is cyclotomic but 
$\chi_2=\psi_K\cdot \chi_{K,2}$ is not.
This implies that $\chi_2$ is non-trivial on $\ker\varphi_2$, so
we have $E_2=0$ and $E(r)=1$.\qed
\medskip\noindent
The preceding proof is remarkably simple in comparison with
the original derivation of (1.4) from (1.3) in [\HO].
The next two sections show that this character sum analysis
generalizes well to more complicated settings.

\head 5. Primes in arithmetic progressions with prescribed primitive root
\endhead

\noindent
For a non-square rational number $r$ as in Theorem 4.2, which is
(under GRH) a primitive root modulo the primes $q$ in a set of positive
density, we now ask what this density becomes
if we restrict to primes $q$ that lie in a prescribed arithmetic
progression.
Thus, on input of $r$, a positive integer $f$ and an integer $a$
coprime to $f$, we want to prove the analogue of Theorem 4.2 for
the set $S=S(r, a\mod f)\subset A$ corresponding to collection of primes
$$
\{
\hbox{$q$ prime : $q\congr a \mod f$ and $r$ is a primitive root modulo $q$}
\}.
\tag{5.1}
$$
The additional congruence condition $q\congr a \mod f$ 
is a condition on the Frobenius of~$q$ in the cyclotomic field
$\Que(\zeta_f)\subset F_\infty=\Que(R_\infty)$. 
In order to formulate it `primewise' at primes dividing~$f$, we use
the natural maps
$$
\rho_p: A_p=\Aut_{R_{p^\infty}\cap\Que^*}(R_{p^\infty})\to
\Aut(\mu_{p^\infty})=\Zee_p^* \to (\Zee_p/f\Zee_p)^*
$$
and take the congruence condition into account by replacing
the primitive root set $A_p\setminus \ker\varphi_p$ at $p$ from (3.2)
by its intersection
$$
S_p = (A_p\setminus \ker\varphi_p) \cap \rho_p^{-1}(a\mod f\Zee_p)
\tag{5.2}
$$
with the congruence set
$\rho_p^{-1}(\overline a)=\rho_p^{-1}(a\mod f\Zee_p)$.
In other words, we map $A_p$ to a finite group by the homomorphism
$$
\varphi_p\times \rho_p: A_p \to
     \Aut_{R_p\cap\Que^*}(R_p)\times (\Zee_p/f\Zee_p)^*,
\tag{5.3}
$$
and let $S_p$ be the inverse image 
$$
S_p =
(\varphi_p\times \rho_p)^{-1}
\left[
(\Aut_{R_p\cap\Que^*}(R_p)\setminus\{1\}) \times \{a\mod f\Zee_p\}
\right].
$$
Note that we have $S_p=A_p\setminus \ker\varphi_p$ 
at primes $p\nmid f$.
With $S_p$ defined as in (5.2), a prime $q>f$ that is coprime to
the numerator and denominator of $r$ and for which 
$\text{Frob}_q\in\Gal(\Que(R_{p^\infty})/\Que)\subset A_p$
lies in $S_p$ for all primes $p<q$
will have $r$ as a primitive root {\it and\/} lie in the
residue class $a \mod f$.

For the rest of this section, we suppose that we are given
coprime integers $a, f\in\Zee_{\ge1}$ and a {\it non-square\/}
rational number $r$, and that $S_p$ is as defined in (5.2).
In this way, $S=\prod_p S_p\subset A$ corresponds to
the set of primes $q$ in (5.1).
\proclaim{5.4. Lemma}
Let $S$ be as defined above, and put
$$
A(r, a\mod f)=
{1\over\phi(f)}
\prod_{p|\gcd(a-1, f)}(1-{1\over p})\cdot
\prod_{p\nmid f} (1-{1\over [F_p:\Que]}),
$$
with $F_p=\Que(\zeta_p, \root p\of r )$ as in $(1.1)$, and
$\phi$ the Euler $\phi$-function.
\item{$1.$}
Suppose $-r$ is not a square.
Then $S$ is non-empty if and only if $r$ is not a $p$-th power for any
prime $p$ dividing $\gcd(a-1, f)$.
In the non-empty case, its measure equals $\nu(S)=A(r, a\mod f)$.
\item{$2.$}
Suppose $-r$ is a square.
Then $S$ is non-empty if and only if the two conditions
\itemitem{\rm (i)}
$r$ is not a $p$-th power for any prime $p$ dividing $\gcd(a-1, f)$;
\itemitem{\rm (ii)}
$a\congr 3\mod 4$ in case $4$ divides $f$;
\item{}
are satisfied.
In the non-empty case, its measure equals $\nu(S)= 2 A(r, a\mod f)$
if\/ $4$ divides $f$, and $\nu(S)= A(r, a\mod f)$ otherwise.
\endproclaim\noindent
\noindent
{\bf Proof.}
We have $\nu(S)=\prod_p \nu_p(S_p)$, and at primes $p\nmid f$
the set $S_p$ has positive measure $\nu_p(S_p)=1-[F_p:\Que]^{-1}$.

Suppose $p$ is an odd prime dividing $f$. 
Then the map $\varphi_p\times\rho_p$ in (5.3) is not surjective, as
it maps $A_p$ onto the fibred product of 
$\Aut_{R_p\cap\Que^*}(R_p)$ and $(\Zee_p/f\Zee_p)^*$ over their
common quotient $\Aut(\mu_p)=(\Zee_p/p\Zee_p)^*$.
For the measure of $S_p$, we have to distinguish two cases.

For $a\not\congr 1\mod p$, the subset $\rho_p^{-1}(\overline a)$ of $A_p$ is
disjoint from $\ker\varphi_p$, so the congruence condition at $p$
implies the primitive root condition at~$p$, and
$S_p=\rho_p^{-1}(\overline a)$ has measure 
$$
\nu_p(S_p)=\nu_p(\rho_p^{-1}(\overline a))
   =(\#(\Zee_p/f\Zee_p)^*)^{-1}=\phi(f_p)^{-1}.
\tag{5.5}
$$
Here $f_p=p^{\ord_p(f)}$ denotes the $p$-part of~$f$.

For $a\congr 1\mod p$, all elements in $\rho_p^{-1}(\overline a)$
are the identity on~$\mu_p$, and we have two subcases.
If $r\in\Que^*$ is a $p$-th~power,
then $\rho_p^{-1}(\overline a)$ is contained in $\ker\varphi_p$ as
we have $\Aut_{R_p\cap\Que^*}(R_p)=\Aut(\mu_p)$.
In this case, $S_p$ and therefore $S$ are empty.
In the other case, in which $r$ is not a $p$-th power in $\Que^*$,
the natural map $\Aut_{R_p\cap\Que^*}(R_p)\to \Aut(\mu_p)$
is $p$ to 1, and we find
$$
\nu_p(S_p)= 
(1-{1\over p})\nu_p(\rho_p^{-1}(\overline a))=
(1-{1\over p}) \cdot \phi(f_p)^{-1}.
\tag{5.6}
$$

For $p=2$ dividing $f$, the map $\varphi_2\times\rho_2$
is surjective in the case where $-r$ is not a square. 
This is because $\Aut(\mu_p)=(\Zee_p/p\Zee_p)^*$ is trivial for $p=2$,
and the action of elements of $A_2$ on $\sqrt r\in R_2$ is
in this case independent of their action on roots of unity of 2-power order.
As $a$ is now odd, we have $2|\gcd(a-1,f)$, and $\nu_2(S_2)$ is given by 
(5.6) for $p=2$ as $r$ is not a square.
This yields the non-twisted case 5.4.1, with $\nu(S)=\prod_p\nu_p(S_p)$
equal to $A(r, a\mod f)$.

In the twisted case 5.4.2 where $-r$ is a square, the action of 
$\alpha\in A_2$ on $\sqrt r\in \i\cdot \Que^*$ and on $\zeta_4=\i$
is `the same', in the sense that we have
$\alpha(\sqrt r)/\sqrt r = \alpha(\zeta_4)/\zeta_4$.
As all $\alpha\in S_2\subset A_2\setminus\ker\varphi_2$ satisfy
$\alpha(\sqrt r)/\sqrt r = -1=\alpha(\zeta_4)/\zeta_4$,
we find that, apart from the necessary condition in (i) for $S$ to
be non-empty, there is the second condition (ii) in case 4 divides $f$.
In the case where we have $4|f$ and $a\congr 3\mod 4$, 
the congruence condition at $2$ implies the primitive root condition at~$2$,
and $\nu_2(S_2)$ is given by (5.5) instead of (5.6) for $p=2$.
Only in this special case, we obtain $\nu(S)=2A(r, a\mod f)$ instead of the
`ordinary' value $\nu(S)=A(r, a\mod f)$.\qed
\medskip\noindent
With the computation of the `naive' density $\nu(S)$ taken care of by 
Lemma 5.4, we can apply Theorem 3.3 to find the actual density
$\delta(S)=\nu (G\cap S)/ \nu(G)$ for the Galois group $G\subset A$ from
Theorem 2.9.
The resulting computation is of striking simplicity when
compared to its original derivation by the second author [\MOa, \MOb]
from a formula analogous to (1.3).
Under GRH, the fraction $\delta(S)$ obtained equals the density,
inside the set of all primes,
of the set of primes $q\congr a\mod f$ for which $r$ is a primitive root.
\proclaim{5.7. Theorem}
Let $a, f\in\Zee_{\ge1}$ be coprime integers,
$r\ne-1$ a non-square rational number that is not a $p$-th power for
any prime $p|\gcd(a-1, f)$.
Define the subset
$S=\prod_p S_p\subset A$ associated to the set of primes in the
residue class $a\mod f$ for which $r$ is a primitive root as in\/ $(5.2)$.
Then we have
$$
\delta(S)={\nu (G\cap S)\over \nu(G)} =
E \cdot
A(r, a\mod f)
$$
for the Galois group $G\subset A$ from $(2.2)$,
with $A(r, a\mod f)$ the Artin constant from Lemma $5.4$,
and the correction factor $E$ equal to
$$
E=1 + E_2\cdot \prod_{p|\gcd(D,f)\text{ odd}} {a\overwithdelims()p}
               \cdot \prod_{p|D,\ p\nmid 2f} {-1\over [F_p:\Que]-1}.
$$
Here $D$ denotes the discriminant of $F_2=\Que(\sqrt r)$, and
$E_2$ is given by
$$
E_2=\cases
-\chi_{F_2, 2}(a)&\text{if $\ord_2(D)\le \ord_2(f)$};\cr
0               &\text{otherwise}.\cr
\endcases
$$
\endproclaim\noindent
%\vglue-.3cm
\noindent
{\bf Proof.}
Suppose first that we are {\it not\/} in the twisted case $D=-4$.
Take $r=r_0^e$ with $e$ odd, so the fields
$K=\Que(\sqrt{r_0})$ from Theorem 2.9
and $F_2=\Que(\sqrt r)$ coincide.
By the assumption on $r$, the naive density
$\nu(S)/\nu(A)$ equals the constant $A(r, a\mod f)$ from Lemma 5.4, 
and we can apply Theorem 3.3 for our set $S$ and the characters 
$\chi_p$ from (3.5) to obtain the correction factor $E=1+\prod_p E_p$.
This amounts to a local computation of $E_p$ at each of the critical
primes~$p|2D$.

At primes $p\nmid 2f$ dividing $D$,
the factors $E_p=-1/([F_p:\Que]-1)$ coming from
the Legendre symbol $\chi_p$ at $p$ are the same
as for Artin's conjecture in Theorem~4.2.
 
For the odd primes $p|\gcd(D,f)$, the Legendre symbol $\chi_p$
has constant value $\chi_p(a)={a\overwithdelims()p}$ on
the congruence set $\rho_p^{-1}(\overline a)$, and therefore on $S_p$.
This yields $E_p={a\overwithdelims()p}$.

Finally, for $p=2$, the character
$
\chi_2=\psi_K\cdot \chi_{K, 2} = \varphi_2 \cdot \chi_{F_2, 2}
$
equals $-\chi_{F_2, 2}$ on $S_2\subset A_2\setminus\ker\varphi_2$.
In the case $\ord_2(D)\le\ord_2(f)$ it has constant value 
$-\chi_{F_2, 2}(a)$ on $S_2\subset \rho_2^{-1}(\overline a)$,
and we obtain $E_2=-\chi_{F_2, 2}(a)$.
In the case $\ord_2(D)>\ord_2(f)$ the character $\chi_{F_2, 2}$
is non-trivial on the subgroup
$\ker\varphi_2\cap \ker \rho_2\subset A_2$.
As $S_2$ is a finite union of cosets of this subgroup, we have
$\int_{S_2} (-\chi_{F_2, 2})d\nu_2=0$, and $E_2=0$.
This finishes the proof in the non-twisted case.

In the twisted case $D=-4$, the field
$K=\Que(\sqrt {r_0})$ from Theorem 2.9 is different from $F_2=\Que(\i)$,
and the correction factor in our theorem simply reads $E=1+E_2$.

If $f$ is not divisible by 4, then $\rho_2$ is the trivial map,
and we have $S_2=A_2\setminus\ker\varphi_2$ and $E_2=0$ as in Theorem 4.2.
In this case we find $\delta(S)=A(r, a\mod f)$.

If 4 divides $f$, we have $1+E_2=1-\chi_{F_2, 2}(a)=1-\chi_{-4}(a)$.
For $a\congr 1\mod 4$ this factor vanishes, and we find 
$\delta(S)=0$, in accordance with the fact that $S$ is empty by~5.4.2.
For $a \congr 3\mod 4$, there is no entanglement correction as
$\chi_2=\psi_K\cdot \chi_{K, 2}$ is non-trivial on the subgroup
$\ker\varphi_2\cap \ker \rho_2\subset A_2$.
We therefore have $\delta(S)=\nu(S)=2A(r, a\mod f)$ by 5.4.2,
and the factor 2 is exactly what $E=1+E_2=2$ gives us.
Note however that in this case, $E$ is a correction for obtaining the right
value of $\nu(S)$, not an entanglement correction factor.\qed
\medskip\noindent
As the Artin constant $A(r, a \mod f)$ is non-zero,
vanishing of the density $\delta(S)$
in Theorem 5.7 occurs if and only if the correction factor $E$ vanishes,
and $G\cap S$ is empty.
It is easy to see when this happens.
\proclaim{5.8. Theorem}
The correction factor $E$ in Theorem $5.7$ vanishes if and only
if we are in one of the two following cases:
\itemitem{{\rm (a)}}
the discriminant of $F_2=\Que(\sqrt r)$ divides $f$, and 
we have $\chi_{F_2}(a)=1$;
\itemitem{{\rm (b)}}
$r$ is a cube in $\Que^*$, the discriminant of
$\Que(\sqrt r)$ divides $3f$ but not $f$, and for $L=\Que(\sqrt{-3r})$
we have $\chi_L(a)=-1$.
\endproclaim\noindent
{\bf Proof.}
The factor $E$ in Theorem 5.7 does not vanish if 
there is a prime $p>3$ that divides the discriminant $D$
of $F_2=\Que(\sqrt r)$ but not $f$,
since then we have $[F_p:\Que]-1\ge p-2>1$.
This leaves us with two cases in which it can vanish.

The first case arises when all odd primes in $D$ divide $f$, and
we have an equality
$
E=1+E_2\prod_{p|D\ \text{ odd}} {a\legendre p}=0.
$
In this case $E_2$ equals $-\chi_{F_2,2}(a)$, so actually $D$ divides~$f$,
and we have $E=1-\chi_{F_2}(a)=0$.
This is case (a), in which all primes congruent to $a\mod f$ are split
in $\Que(\sqrt r)$, making $r$ a quadratic residue modulo 
all but finitely many of these primes.

The second case arises if all odd primes in $D$ divide $f$ except 
the prime $p=3$, which divides $D$ but not $f$,
and we have
$$
E=1+E_2\cdot\prod_{p|D/3\ \text{odd}} {a\legendre p}\cdot
{-1\over [F_3:\Que]-1}=0.
$$
In this situation $E_2$ equals $-\chi_{F_2,2}(a)$,
so $D$ divides $3f$ but not $f$,
and $[F_3:\Que]$ equals~2, showing that $r$ is a cube.
The resulting equality is $E=1+\chi_{L}(a)=0$, so we are in case (b).
To understand this case, we note that a cube can only be a primitive root
modulo a prime $q\congr 2\mod 3$, and no prime $q$
can be inert in all three quadratic subfields of
$\Que(\sqrt r, \sqrt{-3r})$.\qed
\medskip\noindent
The vanishing result 5.8
already occurs in [\LEa, Theorem 8.3], where it is said
to follow from a `straightforward analysis' in terms of Galois groups
that is not further specified.

\head 6. Near-primitive roots
\endhead

\noindent
In addition to $r\in\Que^*\setminus\{\pm1\}$, we 
now let $t=\prod_p t_p\in\Zee_{\ge1}$ be a positive integer,
with $t_p=p^{\ord_p(t)}$ the $p$-component of $p$.
We are interested in the density of the set of primes $q$
for which $r$ is a `near-primitive root' in the sense that $r\mod q$
generates a subgroup of $\FF_q^*$ of exact index $t$.
For odd primes $q$ coprime to numerator and denominator of $r$,
the condition amounts to requiring that
$q$ split completely in the splitting field $F_t=\Que(R_t)$ of $X^t-r$,
but not in any of the fields $F_{pt}$ with $p$ prime.
Note that such primes $q$ will be necessarily be
congruent to $1\mod t$.

In order to define the subset $S=\prod_p S_p\subset A$
for near-primitive roots of index $t$,
we use the surjective restriction maps
$$
\varphi_P: A_p \tto A(P) = \Aut_{R_P\cap\Que^*}(R_P)
\tag{3.1}
$$
for $p$-powers $P$ as defined in Section 3, and put
$$
S_p=\ker \varphi_{t_p} \setminus \ker\varphi_{pt_p}.
\tag{6.1}
$$
Note that (6.1) reduces to (3.2) for $p\nmid t$, when we have $t_p=1$.
Just as in Section~4, $S_p$ is the set-theoretic difference of a group 
and a subgroup.

\proclaim{6.2. Lemma}
Let $P$ be a prime power, and write $r=\pm r_0^e$ as in $(2.3)$.
Then $A(P)$ has order $\phi(P)\cdot P/(P,e)$, unless we are
in the twisted case $r=-r_0^e$ with $P>1$ a $2$-power dividing $e$, when
the order is $2\cdot\phi(P)$.
\endproclaim\noindent
{\bf Proof.}
We can describe the finite quotient $A(P)$ of $A_p$
just as we described its infinite counterpart
$A=\Aut_{R_\infty\cap\Que^*}(R_\infty)$ in Section~2.
Let $r^{1/P}\in\overline\Que$ be any root of the polynomial $X^P-r$. 
Then $R_P=\langle r^{1/P}\rangle \times\mu_P$ is the product
of an infinite cyclic group and the finite group $\mu_P$ of
$P$-th roots of unity, and its quotient
$$
C_P={R_P\over \mu_P\cdot (R_P\cap\Que^*)}
$$
is a finite cyclic group of order dividing $P$,
generated by $r^{1/P} \mod \mu_P\cdot (R_P\cap\Que^*)$.
If $P$ is not a 2-power, or we are not in the twisted case in which
$-r$ is a square, then $r$ is a $(P,e)$-th power in $\Que^*$, and
$C_P$ is of order $P/(P,e)$.
If however $-r=r_0^e$ is a rational square and $P>1$ a 2-power,
$r^{1/P}$ is equal to a primitive $2P$-th root of unity times $r_0^{e/P}$,
and $C_P$ has order $P/(P,e)$ only when $P$ does not divide $e$.
If it does divide $e$, the order is 2 and not $P/(P,e)=1$.

Just as for $A$ and $A_p$, we have an exact sequence
$$
1\to\Hom(C_P,\mu_P)\to A(P)\to \Aut(\mu_P)\to 1
$$
showing that $A(P)$ has order $\phi(P)\cdot \#C_P$.
The result follows.\qed
\medskip\noindent
Using Lemma 6.2, it is straightforward to find the measure of $S_p$ in (6.1),
and the naive density for near-primitive roots.
\proclaim{6.3. Lemma}
Write $r=\pm r_0^e$ as in $(2.3)$.
Then the measure of the set $S=\prod_p S_p$ defined by $(6.1)$ is
equal to
$$
A(r, t) = \alpha_2\cdot{(t,e)\over t^2}\cdot
\prod_{p|t\atop \ord_p(e)\le \ord_p(t)} (1+{1\over p})\cdot
\prod_{p\nmid t}(1-{1\over [F_p:\Que]}),
$$
where $\alpha_2$ is defined by
$$
\alpha_2=\cases
1/2             &\text{if $-r$ is a square and $0<\ord_2(t)\le \ord_2(e)-1$};\cr
1/3             &\text{if $-r$ is a square and $0<\ord_2(t)=\ord_2(e)$};\cr
1               &\text{otherwise}.\cr
\endcases
$$
\endproclaim\noindent
{\bf Proof.}
For primes $p$ that do not divide $t$, we have the familiar Artin factors
$\nu_p(S_p)$.
For primes $p$ dividing $t$, the factors when we are not in the twisted case 
with $p=2$ become
$$
\nu_p(S_p)=
{(t_p,e)\over \phi(t_p)\cdot t_p} -
{(pt_p,e)\over \phi(pt_p)\cdot pt_p} =
{(t_p,e)\over t_p^2}\cdot
\cases
(1 + {1\over p})  &\text{if $\ord_p(e)\le \ord_p(t)$};\cr
1                 &\text{otherwise}.\cr
\endcases
$$
If $2$ divides $t$ and we are in the twisted
case where $-r$ is a square,
we need to correct the value for $\nu_2(S_2)$ given by the formula
above in view of Lemma 6.2.
A short computation shows that this leads to an extra
factor $1/2$ if $2t_2$ divides $e$, and
to a factor $1/3$ if $t_2$ but not $2t_2$ divides $e$.
This is the factor $\alpha_2$.
Taking the product of $\nu_p(S_p)$ over all $p$, we obtain the 
Artin constant $A(r,t)$.\qed

\proclaim{6.4. Theorem}
For $r\in\Que^*$ and $t=\prod_p t_p\in\Zee_{\ge1}$, 
define $S=\prod_p S_p\subset A$ associated to the set of primes 
for which $r$ is a near-primitive root of index $t$ as in\/ $(6.1)$.
Then we have
$$
\delta(S)={\nu (G\cap S)\over \nu(G)} =
E \cdot
A(r, t)
$$
for the Galois group $G\subset A$ from $(2.2)$, with
$A(r, t)$ as in Lemma\/ $6.3$ and
%the correction factor
$E$ equal to
$$
E=1 + E_2\cdot \prod_{p|\disc(K), p\nmid 2t} {-1\over [F_p:\Que]-1}.
$$
Here we take $K=\Que(\sqrt {r_0})$ with $r=\pm r_0^e$ as in $(2.3)$, and
choose $r=r_0^e$ if $e$ is odd.
In terms of $e_2=2^{\ord_2(e)}$ and $d_2=2^{\ord_2(\disc(K))}$,
we have a quantity
$$
s_2=\cases
\lcm(2e_2, d_2)    &\text{if $r=r_0^e$};\cr
4                  &\text{if $-r$ is a square and $(e_2, d_2)=(2,8)$};\cr
4e_2               &\text{if $-r$ is a square and $(e_2, d_2)\ne(2,8)$}\cr
\endcases
$$
that determines the value $E_2$ by
$$
E_2=\cases
1             &\text{if\/ $s_2|t_2$};\cr
0             &\text{if\/ $s_2\nmid 2t_2$};\cr
-1            &\text{if\/ $s_2=2t_2 = 2$};\cr
-1            &\text{if\/ $s_2=2t_2 = 4$, $-r$ is a square
                          and $(e_2, d_2)=(2,8)$};\cr
-1/3          &\text{otherwise}.\cr
\endcases
$$
\endproclaim\noindent
{\bf Proof.}
We already computed $\nu(S)=A(r,t)$ in Lemma 6.3,
so by Theorem 3.3 we only need to check
that $E=1+\prod_{p|2d} E_p$ has the indicated form,
with $d=\disc(K)$.
Note that, even though the field $K=\Que(\sqrt {r_0})$
is only defined up to twisting (as in 2.12.1) by the cyclotomic
character $\chi_{-4}$, divisibility of $d$ by odd primes $p$ 
or by 8 are well-defined notions.
Also, if $e$ is odd, the equality $r=r_0^e$ does uniquely
determine $r_0$ and $K$.

At primes $p|d$ that do not divide $2t$, the factors
$E_p=-1/([F_p:\Que]-1)$ are the same as in 4.1.
At odd primes $p|\disc(K)$ that do divide $t$, the Legendre symbol $\chi_p$
equals 1 on $\ker\varphi_{t_p}$, and therefore on $S_p$.
This yields $E_p=1$ for these $p$, and we obtain the desired expression
$$
E= 1+ E_2 \cdot \prod_{p|\disc(K),\ p\nmid 2t} {-1\over [F_p:\Que]-1},
$$
with $E_2$ the average value of the character
$\chi_2=\psi_K\cdot \chi_{K, 2}: A_2\to\mu_2$ on
the `difference of subgroups'
$S_2=\ker\varphi_{t_2}\setminus\ker\varphi_{2t_2}\subset A_2$.
In order to explicitly find $E_2$, we first compute the smallest 2-power $s_2$
for which $\chi_2$ is trivial on $\ker\varphi_{s_2}$.

By the definitions (2.8) and (3.1) of $\psi_K$ and $\varphi_{2^k}$,
their kernels are the subgroups of $A_2$ that pointwise stabilize 
$\langle r_0^{1/2}\rangle$ and $R_{2^k}$, respectively.
It follows that $\psi_K$ is trivial on $\ker \varphi_{2^k}$
if and only if the group $R_{2^k}\subset\overline\Que^*$
of $2^k$-th roots of $\langle r\rangle$
contains an {\it odd\/} power of $r_0^{1/2}$.
If we are in the untwisted case $r=r_0^e$, we have
$$
R_{2^k}= \langle r_0^{e/2^k}\rangle \times \mu_{2^k} 
\leqno(6.5)
$$
and the smallest 2-power for which this happens is $2^k=2e_2$.
In the twisted case in which $-r$ is a square, we have
$$
R_{2^k}=\langle \zeta_{2^{k+1}}r_0^{e/2^k}\rangle \times  \mu_{2^k}
\leqno(6.6)
$$
for a primitive $2^{k+1}$-th root of unity $\zeta_{2^{k+1}}$, and this
smallest 2-power is $2^k=4e_2$.

For a cyclotomic character of 2-power conductor on $A_2$
such as $\chi_{K, 2}$, it is clear that it is trivial on 
$\ker \varphi_{2^k}$ if and only if its conductor divides~$2^k$.

For $\alpha\in A_2$, the values $\psi_K(\alpha)$ and $\chi_{K, 2}(\alpha)$
depend on the action of $\alpha$ on
$r_0^{1/2}$ and on the $d_2$-th roots of unity.
Thus, $\psi_K$ and $\chi_{K, 2}$ respectively factor via the
`Tate-module' $\Zee_2$ and the cyclotomic component $\Zee_2^*$ of
$A_2\iso\Zee_2\rtimes\Zee_2^*$ (cf. (2.7)).
If we are in the untwisted case (6.5), then 
$\chi_2=\psi_K\cdot \chi_{K, 2}$ is trivial on the
pointwise stabilizer $\ker \varphi_{2^k}$ of $R_{2^k}$
if and only if each of $\psi_K$ and $\chi_{K, 2}$ is, and
we find $s_2=\lcm(2e_2, d_2)$.

Now suppose we are in the twisted case (6.6).
Then $\psi_K$ is trivial on $\ker\varphi_{4e_2}$, and so is $\chi_{K, 2}$
as $d_2$, a divisor of~8, necessarily divides $4e_2$.
The `non-cyclotomic' character $\psi_K$ is non-trivial on
$\ker\varphi_{2e_2}$,
but as every $\sigma\in \ker\varphi_{2e_2}$ fixes
$\zeta_{4e_2}r_0^{1/2}$, it can be described `in cyclotomic terms'
on $\ker\varphi_{2e_2}$ by the identity
$\psi_K(\sigma)=\zeta_{4e_2}^{\sigma-1}$.
Thus, in the case $e_2=2$ and $d_2=8$,
the quadratic characters $\psi_K$ and $\chi_{K, 2}$ are non-trivial
but {\it identical\/} on
$\ker\varphi_{2e_2}=\ker\varphi_4$, so their product
$\chi_2$ is trivial on it.
Apart from this rather special twisted case in which we have $s_2=2e_2=4$,
the character $\chi_2$ is trivial on
$\ker\varphi_{2^k}$ if and only if $\psi_K$ is, i.e.,
if and only if $s_2=4e_2$ divides $2^k$.

Having computed $s_2$, we can easily find $E_2$.
If $s_2$ divides $t_2$, then $\chi_2$ is trivial on $\ker\varphi_{t_2}$,
hence on $S_2$, and we find $E_2=1$.
If $s_2$ does not divide $2t_2$, then
$\chi_2$ is non-trivial on both $\ker\varphi_{2t_2}$ and
$\ker\varphi_{t_2}$, and we find $E_2=0$.
In the remaining case $s_2=2t_2$ we find, just as for 4.1,
$$
E_2={-\nu_2(\ker\varphi_{2t})\over
          \nu_2(\ker\varphi_{t_2}) - \nu_2(\ker\varphi_{2t_2}) }
    ={-1\over [A(t_2):A(2t_2)]-1}.
$$
In the untwisted case $s_2=2t_2=2$, where both $e$ and $d$
are odd, the index $[A_2:A(2)]$ equals 2, and we find $E_2=-1$.
(This is actually a case that already occurred in the proof of 4.2.)
We also find $E_2=-1$ if we have $s_2=2t_2$
in the special twisted case above, for $s_2=2t_2=2e_2=4$ and $d_2=8$;
indeed, we then have $[A(2):A(4)]=2$ from the order formulas
$\#A(2)=2$ and $\#A(4)=4$ provided by Lemma 6.2.
In the other cases with $s_2=2t_2\ge4$ we have $e_2|t_2$, and in the
twisted cases with $s_2=4e_2$ even $2e_2|t_2$.
The order formulas from Lemma 6.2 then yield $\#A(2t_2)=4\cdot\#A(t_2)$ and
$[A(t_2):A(2t_2)]=4$, hence $E_2=-1/3$.
\qed
\medskip\noindent
If we compare Theorem 6.4 to the result for near-primitive root densities
in [\WA], we see that, despite the careful administration we needed
for the twisted case, both the derivation and the resulting expression
for the density given here are considerably simpler.
In fact, it takes some work to see that the formulas in [\WA], which
express the density as a sum of up to 4 different inclusion-exclusion-sums,
can be reduced to our single formula.
Whereas it is extremely cumbersome to derive the vanishing criteria
from the formulas in [\WA], 
it is straightforward to obtain them from~Theorem 6.4.
In the criteria below, which occur without proof as [\LEa, (8.9)--(8.13)],
we write $d(x)$ for $x\in\Zee$
to denote the discriminant of the number field $\Que(\sqrt x)$.
In particular, $d(x)$ equals 1 if $x$ is a square.
\proclaim{6.7. Theorem}
Let $r=\pm r_1^e$ and $t\in\Zee_{\ge1}$ be as in Theorem $6.4$.
Then the near-primitive root density $E\cdot A(r,t)$ in $6.4$
vanishes if and only if we are in one of the following five cases:
\smallskip
\itemitem{{\rm (a)}}
$t$ is odd, and $d(r)|t$;
\itemitem{{\rm (b)}}
$t\congr 2\mod4$, and $r=-u^2$ with $d(2u)|2t$;
\itemitem{{\rm (c)}}
$r$ is a cube, $3\nmid t$, $-r$ is not a square, $d(-3r_0)|t$, and 
$\ord_2(t)>\ord_2(e)$;
\itemitem{{\rm (d)}}
$r$ is a cube, $3\nmid t$, $-r$ is a square, $d(-3r_0)|t$ and
$\ord_2(t)>\ord_2(e)+1$;
\itemitem{{\rm (e)}}
$r$ is a cube, $3\nmid t$, $-r=u^2$, $8|d(-3u)|2t$.
\endproclaim\noindent
{\bf Proof.}
The naive density $A(r,t)$ from Lemma 6.3 vanishes if and only if $t$ is odd
and $r$ is a square. This is case (a) with $d(r)=1$.

As any local factor $E_p=-1/([F_p:\Que]-1)$ satisfies
$|E_p|\le 1/(p-2)<1$ for $p\ge5$, we see that $E=1+\prod_p E_p$
can only vanish if we have $E_p=1$ for all primes $p\ge5$, i.e.,
if all primes $p\ge 5$ dividing $d$ also divide $t$.
Assume that this is the case.
Then $E$ vanishes if and only if we either
have $E_2=-1=-E_3$ or $E_2=1=-E_3$.

Suppose first that we have $E=0$ with $E_2=-1=-E_3$.
For $s_2=2$ this means that $t$ and $e$ and
$d=d(r_0)=d(r)$ are odd, and that $d$ divides $t$.
This is case (a) with $d(r)\ne 1$.
For $s_2=4$ we have $t_2=2$ and $r=-u^2$,
with $u=r_0^{e/2}$ a non-square rational number for which
$d(u)=d(r_0)=d$ satisfies $8|d|4t$.
As $8|d(u)$ can be written as
$\ord_2(d(2u))\le 2=\ord_2(2t)$, we are in case (b).

Suppose next that we have $E=0$ with $E_2=1=-E_3$.
The condition $E_3=-1$ means that
$r$ is a cube, and that 3 divides $d$ but not $t$.
To have $E_2=1$ as well, $t_2$ needs to be divisible by $s_2$,
and this leads to three cases reflecting 
the three cases in the definition of $s_2$.
In the non-twisted case, $t_2$ has to be divisible by $2e_2$ and~$d_2$,
leading to $\ord_2(t)>\ord_2(e)$ and $d(-3r_0)=-d(r_0)/3|t$.
This is case (c).
The twisted case with $s_2=4e_2$ is case (d), with $\ord_2(t)>\ord_2(e)+1$
reflecting the condition $s_2=4e_2|t_2$.
Finally, we have the twisted case with $s_2=4$.
Here $-r=u^2$ is a square and $e_2=2$, so we have $d=d(u)$
and $d(-3u)=-d(u)/3$.
The conditions $4|t$ and $d_2=8$ may now be combined with the
conditions at the odd primes to yield $8|d(-u/3)|2t$,
and we are in case (e).

The reader may check that $E$ indeed vanishes in each of the cases
(a)--(e), or refer to remark 6.8.2 below instead.
\qed
\medskip\noindent
{\bf 6.8. Remarks.}
1. One may restrict case (e) to values $t\congr 4 \mod 8$, as
$t\congr 0\mod 8$ in case (e) is already covered by case (d).
In doing so, the five cases become mutually exclusive.

2. The computation of the vanishing criteria in Theorem 6.7 is so automatic
that one barely realizes {\it why\/} these are vanishing criteria.
In case (a) the number $r$ is a square modulo almost all primes 
$q\congr 1\mod t$, so it cannot generate a subgroup of odd index $t$
modulo such $q$ for $q>2$.
In case (b),
if $r=-u^2$ generates a subgroup of even index modulo $q$,
then ${r\legendre q}={-1\legendre q}=1$ implies that we have
$q\congr 1\mod 4$, and $r=(\i u)^2\mod q$ for a primitive 4-th root
of unity $\i$ modulo $q$.
For $q\congr 1\mod t$ we easily see that $q$ splits in $\Que(\sqrt{u})$
if and only if we have $q\congr 1\mod 8$, so $\i u$ is a square modulo $q$
and $r$ generates a subgroup modulo $q$ of index divisible by 4.

In the cases (c)--(e), the divisibility of the index of $r$ modulo $q$
by $t$ implies that $-3$ is a square modulo $q$, so we have $q\congr 1\mod 3$,
and the cube $r$ generates a subgroup of index divisible by 3.

\head 7. Generalizations
\endhead

\noindent
The examples in the two preceding sections show that the
character sum approach to the computation of various
primitive root densities gives rise to formulas with a simple basic structure.
Unsurprisingly, more case distinctions become necessary as the complexity
of the input data grows.
In more complicated settings, where a single closed formula running over a page
of case distinctions may not be the most desirable result, the method can also
be seen as an {\it algorithm\/} to find the density in each specific case.
\medskip\noindent
{\it Near-primitive roots for primes in arithmetic progressions.}
As a rather straightforward generalization, one may combine Sections 5 and 6
into a single density computation for the set of primes $q\congr a \mod f$
for which a given rational number $r=\pm r_0^e$ generates
a subgroup of exact index $t$ in $\FF_q^*$.
As such primes $q$ are necessarily congruent to $1\mod t$,
it is natural to assume $t|f$ and $a\congr 1\mod t$.
For primes~$p|f$, the original Artin sets 3.2
are then replaced in the spirit of 5.1 and 6.1 by 
$$
S_p= \left(\ker \varphi_{t_p} \setminus \ker\varphi_{pt_p}\right)
           \cap \rho_p^{-1}(a\mod f\Zee_p).
$$
One can now compute the values $\nu_p(S_p)$ and their somewhat
complicated product $A(r, t, a\mod f)$ over all $p$ as before.
Application of Theorem 3.3 yields the fraction
$\delta(S)=\nu(G\cap S)/\nu(G)$ as a product of $A(r, t, a\mod f)$
and a correction factor $E=1+\prod_p E_p$, where the value of $E_2$
requires a large number of case distinctions.
We leave the details to the reader fond of general closed formulas,
and note that when viewed as an {\it algorithm\/}, the method easily yields 
$\delta(S)$ for any set of input values $t, f\in\Zee$,
$(a \mod f)\in(\Zee/f\Zee)^*$ and $r\in\Que^*$.
\medskip\noindent
{\it Higher rank Artin densities.}
There are generalizations of Theorem 3.3
to variants of Artin's conjecture over $\Que$ for which not the
theorem itself, but the general strategy of the proof applies.
One might for instance want to compute, upon input of
$a, b\in \Que^*$, the density of primes $q$ for which $\FF_q^*$ is generated
by $a$ and $b$, or for which $b$ is in the subgroup of $\FF_q^*$ generated
by $a$.
We assume here that we are in the true 2-variable case where $a$ and $b$
are multiplicatively independent, i.e., 
generate a subgroup of rank 2 in $\Que^*/\{\pm1\}$.

In this case, we are led to study the Galois group $G$ of the
extension $\Que\subset\Que(R_\infty)$ obtained by
adjoining to $\Que$ all radicals of $a$ {\it and\/} all radicals of $b$.
Analogously to (2.2), one then has an injective Galois representation
$
G\to A= \Aut_{R_\infty\cap\Que^*} (R_\infty).
$
The group $A$ is an extension of $\Aut(\mu_\infty)=\widehat\Zee^*$ by a free
$\widehat\Zee$-module of rank 2 that naturally decomposes as a product
$A=\prod_p A_p$ of automorphism groups of groups of $p$-power radicals.
The direct analogue of Theorem 2.9 is that $G\subset A$ is a subgroup
of index 4 that arises as the intersection
of the kernels of {\it two\/} quadratic characters
$\kappa=\psi_K\cdot \chi_K$ and $\widetilde\kappa=\psi_{\widetilde K}\cdot \chi_{\widetilde K}$ 
on $A$ related to distinct quadratic entanglement fields $K$
and $\widetilde K$.

For subsets $S=\prod_p S_p\subset \prod_p A_p=A$,
the analogue of Theorem 3.3 is 
that the quotient $\nu(G\cap S)/\nu(G)$ differs from $\nu(S)/\nu(A)$ by
an entanglement correction factor of the form
$$
1+\prod_p E_{\kappa, p} + \prod_p E_{\widetilde\kappa, p} + \prod_p E_{\kappa\widetilde\kappa, p},
$$
with $E_{\alpha, p}$ denoting, for a character $\alpha=\prod_p \alpha_p$ on
$A=\prod_p A_p$, the average value of~$\alpha_p$ on~$S_p$.
It reflects the fact that in this case,
${1\over 4}(1+\kappa+\widetilde\kappa+\kappa\widetilde\kappa)$ is
the characteristic function of $G$ in $A$.
This leads to much easier proofs of results such as
[\MSa, Theorem~3].

Continuing in the direction of arbitrary rank subgroups,
nothing prevents us from considering properties of
subgroups of $\FF_q^*$ that are generated by $n$
elements $a_1, a_2,\ldots, a_n\in\Que^*$ for any $n\in\Zee_{>0}$.
One may for instance look at those $q$ for which {\it all\/} $a_i$
are primitive roots modulo~$q$, or those $q$ for which the subgroup
$\Gamma=\langle a_1, a_2,\ldots, a_n\rangle \subset \Que^*$ 
maps surjectively to~$\FF_q^*$.
Our methods do generalize to this situation, and lead to (short)
proofs and generalizations of theorems obtained previously by Matthews [\MA]
and Cangelmi and Pappalardi~[\CP].
We refer to [\MSb] for further details.

\medskip\noindent
{\it Maximal radical extension.}
The ultimate structural result on the Galois group $G$ over
$\Que$ of the field obtained by adjoining to $\Que$
the group $R=\root\infty\of{\Que^*}$ of all radicals of
{\it all\/} rational numbers is that $G$ is the subgroup of
$A=\Aut_{R\cap \Que^*} (R)$ that is `cut out' by an explicit
family of quadratic characters.
It consists, for each prime $p$, of a character as in Theorem 2.9
that expresses the fact that $\sqrt p$ equals a `Gauss sum',
a sum of roots of unity, and that elements of $G$ should
preserve this {\it additive\/} relation.
It implies that over $\Que$, the groups of radical
Galois extensions $\Que\subset \Que(W)$ for subgroups $W\subset R$
can be described
as subgroups of the automorphism group $\Aut_{\Que^*\cap W}(W)$
that arise as the intersections of kernels of certain quadratic characters.

A beautiful generalization of this result to arbitrary fields $K$
of characteristic zero
was announced in the 2006 lecture notes [\LEb, Section 13]
of the first author.
It describes the Galois group over $K$ of the maximal radical extension 
$K(\root\infty\of {K^*})$ of $K$ explicitly as a subgroup $G$
of $A=\Aut_{K^*}(\root\infty\of {K^*})$.
In all cases, $G$ is normal in $A$, and $A/G$ is a profinite
{\it abelian\/} group.
It opens up the possibility of generalizing all results
that have been proved or mentioned over $\Que$ in this paper
to similar results over arbitrary number fields.
Such extensions, and also generalizations that replace
the multiplicative group by one-dimensional tori,
are the subject of Palenstijn's PhD-thesis [\WJP].

\head 8. Cyclic reduction of Serre curves
\endhead

\noindent
There is no end to the number of generalizations of Artin's conjecture
that one may wish to consider, and we conclude this paper by an
example illustrating the applicability of our character
sum method in an {\it elliptic\/} setting.
To stress the analogy, we will re-use the notation of Section~2
($G$, $A$, $F_\infty$, $\ldots$) in this final section to denote 
elliptic analogues of the earlier objects.

Let $E$ be an elliptic curve defined over $\Que$,
and $F_\infty=\Que(E^\tor(\overline\Que))$ the infinite extension of~$\Que$
generated by the coordinates of the torsion points of $E$.
As the absolute Galois group of $\Que$ acts by group automorphisms on
$E^\tor(\overline\Que)$, we have an injective Galois representation
$$
G=\Gal(F_\infty/\Que) \longrightarrow A=\Aut(E^\tor(\overline\Que)).
\tag{8.1}
$$
Just as (2.2) is obtained from finite group homomorphisms (2.1) for
$n\in\Zee_{\ge1}$,
the map (8.1) is obtained as the profinite limit of finite
Galois representations
$$
\Gal(F_n/\Que)\tto
A(n)=\Aut (E[n](\overline\Que)) \iso\GL_2(\Zee/n\Zee)
\tag{8.2}
$$
describing the Galois group of the $n$-division field
$F_n=\Que(E[n](\overline\Que))$ of $E$ for $n\in\Zee_{\ge1}$.
Picking compatible isomorphisms $A(n)\iso\GL_2(\Zee/n\Zee)$,
we obtain an isomorphism of profinite groups
$A\isar \GL_2(\widehat\Zee)$ and
a map
$$
\omega: A \isar\GL_2(\widehat\Zee) \mapright{\det}\widehat\Zee^*
$$
known as the {\it cyclotomic character\/}.
The field $F_\infty$ contains $\Que(\mu_\infty)=\Que_\ab$, 
and the action of $\sigma\in G$
on roots of unity $\zeta\in\mu_\infty$ is given by
$\sigma(\zeta)=\zeta^{\omega(\sigma)}$.

Serre [\SE] proved that for $E$ without CM,
the {\it image of Galois\/} $G$ is an open subgroup of~$A$,
and hence of finite index in $A$.
He also observed that we always have $G\ne A$, as there is a
non-trivial quadratic character $\chi_E: A\to\mu_2$ associated to
$E$ that vanishes on $G$.

In order to describe Serre's character, let
$\psi: A\to A(2)\to \mu_2$
be the non-trivial quadratic character that maps $\alpha\in A$
to the sign of the permutation by which it acts on the three 
points of order 2 on $E$.
If we write $E$ by an explicit Weierstrass model $y^2=f(x)$
over $\Que$, with $f\in\Que[x]$ cubic of discriminant $\Delta\ne0$,
these three points of order 2 have coordinates $(e,0)$,
with $e\in\overline\Que$ a zero of $f$.
The 2-division field $F_2$ is the splitting field of $f$,
and for a field automorphism $\sigma\in G$,
the sign of the permutation action of $\sigma$ on the points of order 2 
is equal to the action
$$
\chi_K(\sigma)=\sigma(\sqrt \Delta)/\sqrt \Delta)\in\mu_2
$$
of $\sigma$ on the subfield $K=\Que(\sqrt \Delta)\subset F_2$.
Note that $K$ does not depend on the Weierstrass model for $E$.
We have $K=\Que$ in the case where the 2-division field $F_2$
is a cyclic cubic field, or equal to $\Que$ itself.
We extend $\chi_K$ to the quadratic character $\chi_K: A\to\mu_2$
on $A$ that is obtained by composing the cyclotomic character
$\omega: A\to \widehat\Zee^*$ with the Dirichlet character associated to $K$.
Then $\chi_K$ is trivial exactly when $\Delta$ is a square,
and on the image of Galois, the characters $\psi$ and $\chi_K$ coincide.
Thus $G$ is contained in the kernel of the quadratic character
$\chi_E=\psi\cdot\chi_K$, and $\chi_E$ is non-trivial as
$\chi_K$ factors via the cyclotomic character $\omega$, whereas
$\psi$ does not.
Just as in (3.5), the character $\chi_E$ can be written as a product 
of characters $\prod_p\chi_p$ on $A=\prod_p A_p$,
with $A_p=\Aut(E[p^\infty](\overline\Que))$.
Here $\chi_p: A\mapright{\omega}\widehat\Zee^*\to\mu_2$ is the Legendre
symbol at $p$ for the odd primes $p$ dividing $D=\disc(K)$,
and $\chi_2=\psi\cdot\chi_{K,2}$ is a non-trivial 
quadratic character that factors via $A(2)$ if and only if 
$D$ is odd.

A {\it Serre curve\/} is an elliptic curve over $\Que$ for
which $G=\ker\chi_E$ is of minimal index $[A:G]=2$.
For Serre curves, the $p$-division field $F_p$ is clearly
of maximal degree
$$
[F_p:\Que]=\#\GL_2(\FF_p)=(p^2-1)(p^2-p)
$$
for all odd primes $p$.
In fact, the same is true for $p=2$, as in this case $F_2$ cannot be
cyclic cubic
(that would make $F_2$ cyclotomic, and $[A:G]$ divisible by 3)
or equal to $\Que$. 
In particular, $\chi_K$ is non-trivial and $D\ne 1$
in the case of Serre curves.

Nathan Jones [\JO] recently proved that, in an asymptotic sense
that is easily made precise, {\it almost all\/}
elliptic curves over $\Que$ are Serre curves.
In this generic case, the image of Galois $G$
inside $A$ is cut out by a single quadratic character
$\chi_E = \psi\cdot\chi_K$, just like in our Theorems 2.9 and 3.3.

For an elliptic curve $E$ over $\Que$ and $q$ a prime of
good reduction, the group $E(\FF_q)$ is a finite group that
is cyclic if and only if for no prime $p<q$, the full
$p$-torsion of $E$ is defined over $\FF_q$.
Equivalently, the prime $q$ does not split completely
in de $p$-division field $F_p$ of $E$.
Serre showed in 1978 that Hooley's argument [\HO] mentioned
in the Introduction can be adapted to show that also in this case,
under GRH, the density of the set of
primes $q$ for which $E(\FF_q)$ is cyclic is
given by the inclusion-exclusion sum
$$
\delta(E)= \sum_{n=1}^\infty {\mu(n)\over [F_n:\Que]}
\tag8.3
$$
that we already encountered for Artin's conjecture as (1.3).
Nobody ever seems to have taken the trouble of explicitly
evaluating such densities, or relating them to the elliptic Artin constant
$$
{\Cal A}_{\text {ell}}=
\prod_{p\text{ prime}} \Bigl(1-{1\over (p^2-1)(p^2-p)}\Bigr)
\doteq .8137519 .
$$
Formula (8.3) holds for any elliptic curve $E$ over $\Que$.
For Serre curves, we have immediate elliptic analogues of Theorems
2.9 and 3.3, and our method yields the following 
analogue of Theorem 4.2.
For lack of a twisted case, the proof is even shorter than it was for 
Artin's conjecture.
\proclaim
{8.4. Theorem}
Let $E$ be a Serre curve, and $D=\disc(\Que(\sqrt\Delta))$ as above.
Then the value $\delta(E)$ in $(8.3)$ is equal to
$$
\delta(E) = C_E\cdot \prod_{p\text{ prime}}
                 \Bigl(1-{1\over [F_p:\Que]}\Bigr)
        = C_E\cdot {\Cal A}_{\text {ell}}
%\prod_{p\text{ prime}} \Bigl(1-{1\over (p^2-1)(p^2-p)}\Bigr)
$$
for an entanglement correction factor $C_E$ that has the value $1$
if $D$ is even, and the value
$$
C_E=1 +\prod_{p|2D} {-1\over [F_p:\Que]-1}
=1 +\prod_{p|2D} {-1\over (p^2-1)(p^2-p)-1}
$$
if $D$ is odd.
\endproclaim
\noindent
{\bf Proof.}
As for 4.2.
When $C_E$ is written in the form (1.5) provided by Theorem 3.3,
it suffices to observe that the local factor $E_2$ is given by (4.1)
if and only if the local factor $\chi_2$ of $\chi_E$ factors via $A(2)$,
i.e., if and only if $D$ is odd. For even $D$ we have $E_2=0$ and $C_E=1$.
\qed
\medskip\noindent
We refer to [\BR] for similar formulas for elliptic curves
that are not assumed to be Serre curves.

\enddocument